\documentclass[a4paper]{article}
\usepackage{caption}
\usepackage{graphicx}
\usepackage{psfrag}
\usepackage{epsfig}
\usepackage[pass]{geometry}
\usepackage{amssymb}
\usepackage{amsmath}
\usepackage{amsfonts}
\usepackage{amssymb}
\usepackage{mathtools}
\usepackage{array}
\usepackage{hyperref}
\usepackage[makeroom]{cancel}
\usepackage{subfig}
\usepackage{float}
\usepackage{calc}
\usepackage{tikz}
\usetikzlibrary{calc}
\usepackage{algorithm}
\usepackage{algorithmic}
\usepackage{placeins}
\newcommand{\R}{\mathbb{R}}
\newcommand{\Z}{\mathbb{Z}}
\newcommand{\N}{\mathbb{N}}
\newcommand{\Lspace}[1]{\mathrm{L}^{#1}}
\newcommand{\C}[1]{\mathrm{C}^{#1}}
\newcommand{\vmin}{v_{\min}}
\newcommand{\vmax}{v_{\max}}
\newcommand{\fmin}{f_{\min}}
\newcommand{\fmax}{f_{\max}}
\newcommand{\rhomin}{\rho^{\min}_0}
\newcommand{\rhomax}{\rho^{\max}_0}
\newcommand{\rhom}{\rho_{\max}}
\newcommand{\rhocr}{\rho_{\mathrm{cr}}}
\newcommand{\Outflow}{\textrm{Out}}
\newcommand{\Inflow}{\textrm{In}}
\newcommand{\ftarget}{f^*}
\newcommand{\vhat}{\hat{v}}
\newcommand{\BV}{\mathrm{BV}}
\newcommand{\TV}{\mathrm{TV}}
\newcommand{\dx}{\Delta x}
\newcommand{\dt}{\Delta t}
\newcommand{\dv}{\Delta v}
\newcommand{\xj}[2]{x_{#1\frac{#2}{2}}}
\newcommand{\modulo}[1]{{\left|#1\right|}}
\newcommand{\norma}[1]{{\left\|#1\right\|}}

\newtheorem{theorem}{Theorem}
\newtheorem{lemma}[theorem]{Lemma}

\newtheorem{proposition}{Proposition}
\newtheorem{definition}{Definition}[section]
\newtheorem{remark}{Remark}

\newtheorem{problem}{Problem}
\newtheorem{hypothesis}{Hypothesis}
\newenvironment{Proof}{{\it Proof.}~}{\hfill$\square$}

\title{Traffic regulation via controlled speed limit%
		\thanks{This research was supported by the NSF grant CNS $\#$1446715 and by KI-Net "Kinetic description of emerging challenges in multiscale problems of natural sciences" - NSF grant \# 1107444}}
\author{Maria Laura Delle Monache%
		\thanks{Department of Mathematical Sciences, Rutgers University - Camden, Camden, NJ, USA (email: \href{mailto:ml.dellemonache@rutgers.edu}{ml.dellemonache@rutgers.edu}).}
		\and
		Benedetto Piccoli%
		\thanks{Department of Mathematical Sciences and CCIB, Rutgers University - Camden, Camden, NJ USA (email: \href{mailto:piccoli@camden.rutgers.edu}{piccoli@camden.rutgers.edu}).}
		\and
		Francesco Rossi%
		\thanks{Aix Marseille Universit\'e, CNRS, ENSAM, Universit\'e de Toulon, LSIS UMR 7296,13397, Marseille, France (email: \href{mailto:francesco.rossi@lsis.org}{francesco.rossi@lsis.org}).}}
\date{}

\begin{document}
\maketitle

\begin{abstract}
We study an optimal control problem for traffic regulation via variable speed limit. The traffic flow dynamics is described with the Lighthill-Whitham-Richards (LWR) model with Newell-Daganzo flux function. We aim at minimizing the $\Lspace{2}$ quadratic error to a desired outflow, given an inflow on a single road. We first provide existence of a minimizer and compute analytically the cost functional 
variations due to needle-like variation in the control policy. Then, we compare three strategies: instantaneous policy; random exploration of control space; steepest descent using numerical expression of gradient. We show that the gradient technique is able to achieve a cost within 10\% of random exploration minimum with better computational performances.
\end{abstract}

\section{Introduction}\label{sec:intro}
In this paper, we study an optimal control problem for traffic flow on a single road using a variable speed limit. The first traffic flow models on a single road of infinite length using a non-linear scalar hyperbolic partial differential equation (PDE) are due to  Lighthill and Whitham \cite{LW55} and, independently, Richards \cite{R56}, which in the 1950s proposed a fluid dynamic model to describe traffic flow. Later on, the model was extended to networks \cite{GPbook}  and started to be used to control and optimize traffic flow on roads. In the last decade, several authors studied optimization and control of conservation laws and several papers proposed different approaches to optimization of hyperbolic PDEs, see \cite{BH12, FHKM06,  GU10, GHKL05, JKCDW06, U02, U03} and references therein. These techniques were then employed to optimize traffic flow through, for example, inflow regulation \cite{CG07}, ramp-metering \cite{RKDMSGB14} and variable speed limit \cite{GGK14}. We focus on the last approach, where the control is given by the maximal speed allowed on the road. Notice that also the engineering literature presents a wealth of approaches \cite{ADFFP98, ADFFP99, CDW11, CPPM10, CVH13, FC12, HDSH05_2, HDSH05, HH10, HHSS09, HHSSV08, HXZ07, YLLZ13}, but mostly in the time discrete setting. In \cite{ADFFP98, ADFFP99} a dynamic feedback control law is employed to compute variable speed limits using a discrete macroscopic model. Instead, \cite{HDSH05_2, HDSH05, HH10} use model predictive control (MPC) to optimally coordinate variable speed limits for freeway traffic with the aim of suppressing shock waves.\\
In this paper, we address the speed limit problem on a single road.
The control variable is the maximal allowed velocity, which may vary in time but we assume
to be of bounded total variation, and
we aim at tracking a given target outgoing flow. More precisely, the main goal is to minimize the quadratic difference between the achieved outflow and the given target outflow. 
Mathematically the problem is very hard, because of the delays in the effect of the control variable (speed limit). In fact, the Link Entering Time (LET) $\tau(t)$, which represents the entering time of the car exiting the road at time $t$ see \eqref{eq:tau},  
depends on the given inflow and the control policy on the whole
time interval $[\tau(t),t]$. Moreover, the input-output map is defined in terms of LET,
thus the achieved outflow at time $t$ depends on the control variable on the whole interval $[\tau(t),t]$.
Due to the complexity of the problem, in this article we restrict the problem to free flow conditions. Notice that this assumption is not too restrictive. Indeed, if the road is initially in free flow, then it will keep the free flow condition due to properties of the LWR model, see \cite[Lemma 1]{BNP06}.\\
After formulating the optimal control problem, we consider needle-like variations for the control policy
as used in the classical Pontryagin Maximum Principle \cite{BPbook}.
We are able to derive an analytical expression of the one-sided variation of the cost,
corresponding to needle-like variations of the control policy, using fine properties
of functions with bounded variation. In particular the one-sided variations depend
on the sign of the control variation and involves integrals w.r.t. to the distributional derivative
of the solution as a measure, see \eqref{eq:gradient}. 
This allows us to prove Lipschitz continuity of the cost
functional in the space of bounded variation function and prove existence of a solution.\\
Afterwards, we define three different techniques to solve numerically this problem.
\begin{itemize}
\item Instantaneous Policy (IP). We design a closed-loop policy, which depends only
on the instantaneous density at road exit. More precisely, we choose the speed limit which gives the nearest outflow to the desired one.
\item Random Exploration (RE). It uses time discretization and random binary tree search of the control space to find the best maximal velocity profile.
\item Gradient Descent Method (GDM). It consists in approximating numerically the gradient of the cost functional using \eqref{eq:gradient} combined with a steepest descent method.
\end{itemize}
We compare the three approaches on two test cases: constant desired outflow and sinusoidal inflow; sinusoidal desired outflow and inflow. In both cases RE provides the best control policy, however GDM performs within 10\% of best RE result with a computational cost of around 15\% of RE. On the other side, IP performs poorly with respect to the RE, but with a very low computational cost. Notice that, in some cases, IP may be the only practical policy, while GDM represents a valid approach also for real-time control, due to good performances and reasonable computational costs. Moreover, control policies provided by RE may have too large total variation to be of practical use.\\
The paper is organized as follows:  section \ref{sec:math_model} gives the description of the traffic flow model and of the  optimal control problem. Moreover, the  existence of
a solution is proved. In section \ref{sec:policies}, the three different approaches to find control policies
are described. Then in section \ref{sec:num_simulations}, these techniques are implemented on two test cases. Final remarks and future work are discussed in section \ref{sec:conclusions}.
\section{Mathematical model}\label{sec:math_model}
 In this section, we introduce a mathematical framework for the speed regulation
 problem. The traffic dynamics is based on the classical Lighthil-Whitham-Richards (LWR) model (\cite{LW55, R56}), while the optimization problem will
 seek minimizers of quadratic distance to an assigned outflow.
\subsection{Traffic flow modeling}\label{subsec:traffic_model}
\begin{figure}[ht]
	\centering
	\subfloat[Velocity function.]{	
	\resizebox{.4\columnwidth}{!}{
	\begin{tikzpicture}[scale=1.4,baseline=2]
\draw [<->] (0,3.5) -- (0,0) -- (5,0);
\draw (4.5,0)parabola(2,2);
\draw (2,2)--(0,2); 
\draw [dashed](2,1.3)--(0,1.3); 
\draw [dashed](4.5,0)parabola(2,1.3);
\draw [dashed](2,2.7)--(0,2.7); 
\draw [dashed](4.5,0)parabola(2,2.7);
\draw [dotted] (2,2.7)--(2,0);
\node [below] at (4.9,0) {$\rho$};
\node [below] at (2,0) {$\rho_{\mathrm{cr}}$};
\node [below] at (4.5,0) {$\rho_{\max}$};
\node [left] at (0,2) {$v(t)$};
\node [left] at (0,1.3) {$v_{\min}$};
\node [left] at (0,2.7) {$v_{\max}$};
\node [left] at (0,3.3) {$v(\rho)$};
\end{tikzpicture}}
	\label{fig:speed}}
	\subfloat[Newell-Daganzo fundamental diagram.]{	
	\resizebox{.45\columnwidth}{!}{
	\begin{tikzpicture}[scale=1.4,baseline=2]
\draw[<->](0,4)--(0,0)--(5,0);
\draw[thick](0,0)--(2,2)--(4.5,0);
\draw[dashed, thick](0,0)--(2,1)--(4.5,0);
\draw[dashed,thick](0,0)--(2,3)--(4.5,0);
\draw[dotted](2,0)--(2,2);
\draw[dotted](0,2)--(2,2);
\draw(1.3,1.3)--(1.3,1)--(1,1);
\node[below]at(5,0){$\rho$};
\node[below]at(2,0){$\rho_{\mathrm{cr}}$};
\node[below]at(4.5,0){$\rho_{\max}$};
\node[left]at(0,4){$f(\rho)$};
\node[left]at(0,2){$f(\rho_{\mathrm{cr}})$};
\node[right]at(1.3,1.2){$v(t)$};
\node[right]at(1.,0.4){$v_{\min}$};
\node[left]at(1.6,2.5){$v_{\max}$};
\end{tikzpicture}}
	\label{fig:FD_triang_VSL}}
	\caption{Velocity and flow for different speed limits.}
	\label{fig:VSL_FD}
\end{figure}
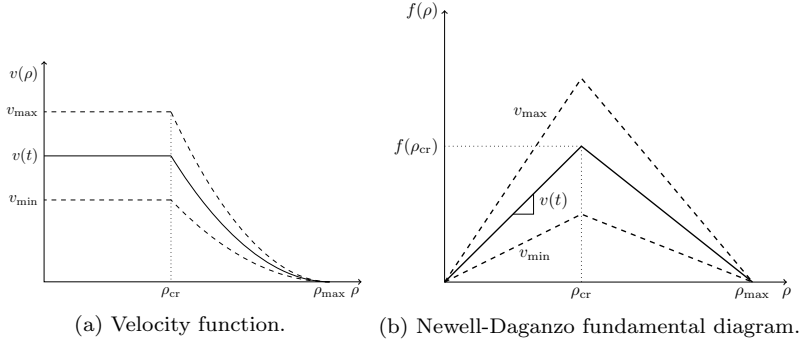

We consider the LWR model  on a single road of length $L$ to describe the traffic dynamics. 
The evolution in time of the car density $\rho$ is described by a 
Cauchy problem for scalar conservation law with time dependent 
maximal speed $v(t)$:
\begin{equation}
	\left\lbrace
	\begin{array}{ll}
	\rho_t + f(\rho, v(t))_x = 0, &\qquad (t,x)\in \R^+\times [0,L],\\
	\rho(0,x) = \rho_0(x), &\qquad x\in [0,L],
	\end{array}
	\right. 
	\label{eq:cons_law}
\end{equation} 
where  $\rho=\rho(t,x)\in [0,\rhom]$ with $\rhom$ the maximal car density. 
In the transportation literature the graph of the flux function $\rho\to f(\rho)$ 
(in our case for a fixed $v(t)$) is commonly referred to as the fundamental diagram. 
Throughout the paper, we focus on the Newell - Daganzo - type (\cite{D94}) fundamental diagrams, see Figure \ref{fig:FD_triang_VSL}.
 The speed takes value on a bounded interval $v(t) \in [\vmin , \vmax],$ $0< \vmin \leq \vmax$, thus the flux function 
$f: [0,\rhom]\times [\vmin , \vmax]\rightarrow \R^+$ is given by 
\begin{equation}
\label{eq:flux}
	f(\rho, v(t))=
	\left\lbrace
	\begin{array}{ll}
	\rho v(t) , & \qquad \text{if } 0 \leq \rho \leq \rhocr, \\
	\dfrac{v(t)\rhocr}{\rhom-\rhocr}(\rhom -\rho), &\qquad \text{if } \rhocr < \rho \leq \rhom , 
	\end{array}
	\right.
\end{equation}
with $v(t)$  representing the maximal speed, see Figure \ref{fig:speed}.
Notice that the flow is increasing up to a \emph{critical density} $\rhocr$
and then decreasing. The interval $[0,\rhocr]$ is referred to as
the \emph{free flow} zone, while $[\rhocr,\rhom]$ is referred to
as the \emph{congested flow} zone.\\

The problem we consider is the following. Given an inflow $\Inflow (t)$,
we want to track a fixed outflow $\Outflow (t)$ on a time horizon $[0,T],$ $T>0$,  by acting on the time-dependent maximal velocity $v(t).$ 
A maximal velocity function $v:[0,T]\rightarrow [\vmin,\vmax]$ 
is called a \textbf{control policy}.\\
It is easy to see that a road in free flow can become congested
only because of the outflow regulation with shocks moving backward,
see \cite[Lemma 2.3]{BNP06}. Since we assume Neumann boundary conditions at the road exit,
the traffic will always remain in free flow, i.e.
$\rho(t,x) \leq \rhocr$ for every $(t,x)\in [0,T]\times [0,L].$ 
Given the inflow function $\Inflow(t)$, we consider the Initial Boundary Value Problem
with assigned flow boundary condition $f_l$ on the left and Neumann boundary condition (flow $f_r$) on the right:
\begin{equation}
	\left\lbrace
	\begin{array}{ll}
	\rho_t + f(\rho, v(t))_x = 0, &\qquad (t,x)\in \R^+\times [0,L],\\
	\rho(0,x) = \rho_0(x), &\qquad x\in [0,L],\\
	 f_{l}(t)=\Inflow(t),\\
	f_r(t)=\rho(t,L)\, v(t).
	\end{array}
	\right. 
	\label{eq:cons_law_IBVP}
\end{equation} 
We denote by $\BV$ the space of scalar functions of bounded variations and by $\TV$ the total variation, see \cite{Bbook} for details. 
For any scalar $\BV$ function $h$ we denote by $\xi(x^\pm)$ its right
(respectively left) limit at $x$.
We further assume the following:
\begin{hypothesis}\label{Hyp:BV}
There exists $0<\rhomin\leq \rhomax\leq\rhocr$ and $0<\fmin \leq \fmax$ such that $\rho_0 \in \BV([0,L],[\rhomin,\rhomax])$ and $\Inflow \in \BV([0,T],[\fmin,\fmax]).$
\end{hypothesis}
Under this assumption, we have: 
\begin{proposition}\label{prop:exi_uni}
Assume that hypothesis \ref{Hyp:BV} holds and $$v \in \BV([0,T],[\vmin,\vmax]).$$
Then, there exists a unique entropy solution $\rho(t,x)$ to \eqref{eq:cons_law_IBVP}.
Moreover, $\rho(t,x)\leq \rhocr$ and,
setting 
\begin{equation}\label{eq:Out}
\Outflow(t)=\rho(t,L)v(t), 
\end{equation}
we have that $\Outflow(.)\in \BV([0,T], \R)$ and
the following estimates hold
\begin{equation}
	\label{eq:estimate_rho}
	\min {\Big\{ \rhomin, \dfrac{\fmin}{\vmax} \Big\}}  \leq \rho (t) \leq \max{\Big\{\rhomax, \dfrac{\fmax}{\vmin} \Big\}},
\end{equation}
\begin{equation}
	\label{eq:estimate_flux}
	\min {\Big\{ \rhomin \vmin, \fmin\dfrac{\vmin}{\vmax} \Big\}}  \leq \Outflow(t) \leq \max{\Big\{\rhomax \vmax, \fmax\dfrac{\vmax}{\vmin} \Big\}}.
\end{equation}
\end{proposition}

\begin{figure}[ht]
\centering
\begin{tikzpicture}[scale=1.4,baseline=0]
\draw[thick](0,0)--(4,0);
\draw[thick](0,0)--(0,4);
\draw[thick](4,0)--(4,4);
\draw[white, fill=gray!20!white, opacity=0.5](0,0)rectangle(4,4);
\draw[line width=1pt] (0,0)--(1.5,0.8)--(3,1.3)--(4,1.8);
\draw[line width=1pt] (0,2)--(1.5,2.5)--(2.5,3)--(4,3.5);
\node[below]at(0,0){$0$};
\node[below]at(4,0){$L$};
\node[left]at(0,2){$\tau(t)$};
\node[right]at(4,1.8){$t_0$};
\node[right]at(4,3.5){$t$};
\end{tikzpicture}
\caption{Graphical representation of the LET function $\tau=\tau(t,v)$ defined in \eqref{eq:tau}.}
\label{fig:tau}
\end{figure}
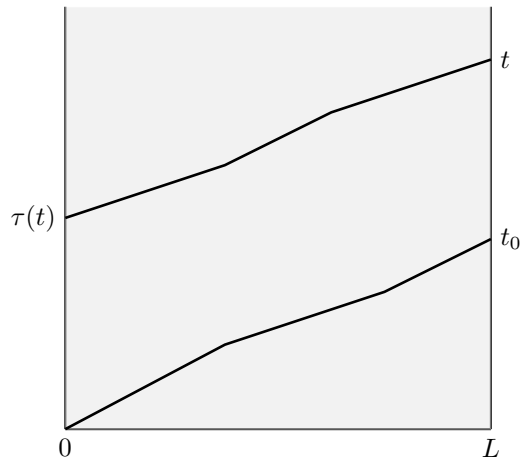

\begin{Proof}
Let $v^n \in \BV([0,T],[\vmin, \vmax])$ be a sequence of piecewise constant functions converging to $v$ in $\Lspace{1}$ and satisfying $\TV(v^n)\leq \TV(v).$ 
For each $v^n$, by standard properties  of Initial Boundary Value Problems for conservation laws \cite{DM07}, 
there exists a unique BV entropy solution $\rho^n$ to \eqref{eq:cons_law_IBVP}
with $\rho^n \in \mathrm{Lip} ([0,T],\Lspace{1})$. 
Notice that the left flow condition is equivalent to
the boundary condition: $\rho_l(t)=\dfrac{\Inflow(t)}{v(t)}$.
From \cite[Lemma 2.3]{BNP06} and the Neumann boundary condition on the right, we get
that $\rho^n(t,x)\leq \rhocr$, thus by maximum principle it holds:
$$
\rho^n(t,\cdot)\in
\BV\Big( \R, \Big[ \min{\Big\{  \rhomin ,\; \dfrac{\fmin}{\vmax} \Big\} }, \max{\Big\{  \rhomax, \; \dfrac{\fmax}{\vmin}\Big\}} \Big] \Big).
$$
Let us now estimate the total variation of the solution $\rho^n$.
Since it solves a scalar conservation laws, the total variation does not
increase in time due to dynamics on $]0,L[$.
Notice that changes in $v(\cdot)$ will not increase the total variation
of $\rho^n$ inside the road (i.e. on $]0,L[$). The total variation of $\rho^n$
increases only because of new waves generated by changes in the inflow.
Using the boundary condition  $\rho_l(t)=\dfrac{\Inflow(t)}{v(t)}$,
we can estimate the total variation in space of $\rho^n$ caused by time
variation of $\Inflow$, respectively time variation of $v$, by $\dfrac{\TV(\Inflow)}{\vmin}$,
respectively $\dfrac{\fmax\,\TV(v)}{\vmin^2}$. Finally we get:
\[
\sup_t \TV (\rho^n(t,\cdot))\leq
\TV(\rho^n(0,\cdot)+\dfrac{\TV(\Inflow)}{\vmin}+
\dfrac{\fmax\,\TV(v)}{\vmin^2}.
\]
By Helly's Theorem (see \cite[Theorem 2.4]{Bbook}) there exists a subsequence converging  in $\Lspace{1}([0,T]\times [0,L])$ to a limit $\rho^*$.
By Lipschitz continuity of the flux and dominated convergence we get that
$f(\rho^n(t,x),v(t))$ converges 
in $\Lspace{1}([0,T]\times [0,L])$ to $f(\rho^*(t,x),v(t))$.
Passing to the limit in the weak formulation
$\int_{\Omega} \rho^n\,\varphi_t+f(\rho^n,w)\,\varphi_x\ dt\,dx=0$
(where $\Omega\subset\subset [0,T]\times [0,L]$ and
$\varphi\in\C{\infty}_0$) we have that $\rho^*$ is a weak entropic solution.
We can pass to the limit also in the left boundary condition because 
this is equivalent to $\rho_l(t)=\dfrac{\Inflow(t)}{v(t)}$ and $v$
is bounded from below.
Finally $\rho^*$ is a solution to \eqref{eq:cons_law_IBVP}. 
 The standard Kru\v{z}hkov entropy condition \cite{K70} ensures uniqueness of the solution.\\
Since $\Outflow(t) = \rho(t,L)v(t)$, we have that $\Outflow(t)$ 
has bounded variation and satisfies  \eqref{eq:estimate_flux}.
\end{Proof}
To simplify notation,
we further make the following assumptions: 
\begin{hypothesis}\label{hyp:bounds}
We assume Hypothesis \ref{Hyp:BV} and the following:
 $$
\rhomin \leq \dfrac{\fmin}{\vmax} \qquad \text{and} \qquad \rhomax \geq \dfrac{\fmax}{\vmin}.
 $$ 
 \end{hypothesis}
 Given a control policy $v$,
we can define a Link Entering Time (LET) function $\tau=\tau(t,v)$ representing the entering time for a car exiting the road at time $t$. The function depends on the control policy $v$, but for simplicity we will write $\tau(t)$ when the policy is clear from the context. Notice that LET is defined only for time greater than a given $t_0>0$, the exit time of the car entering the road at time $t=0$, see Figure \ref{fig:tau}.
Note that $t_0$ satisfies $\int_0^{t_0}v(s)ds=L$ and,
for each $t\geq t_0$:
\begin{equation}
	\label{eq:tau}
	\int_{\tau(t)}^{t} v(s)ds = L.
\end{equation}
Such $\tau(t)$ is unique, due to the hypothesis $v\geq \vmin >0.$
From the identity 
$$
\int_{\tau(t_1)}^{\tau(t_2)}v(s)ds = \int_{t_1}^{t_2}v(s)ds,
$$
we get the following:
\begin{lemma} \label{lem:tau}
Given a control policy $v$, the function $\tau$ is a Lipschitz continuous function, with Lipschitz constant $\dfrac{\vmax}{\vmin}$.
\end{lemma}
Recalling the definition of outflow of the solution given in \eqref{eq:Out}, we get:
\begin{proposition}\label{prop:inputOutputMap}
The input-output flow map of the Initial Boundary Value Problem (briefly IBVP)
\eqref{eq:cons_law_IBVP} is given by 
\begin{equation}
	\label{eq:Outflow}
	\Outflow(t)=\Inflow(\tau(t))\dfrac{v(t)}{v(\tau(t))}.
\end{equation}
\end{proposition}
\begin{Proof}
Thanks to Proposition \ref{prop:exi_uni}, the solution $\rho$ to the IBVP 
\eqref{eq:cons_law_IBVP} satisfies $\rho(t,x)\leq \rhocr$, thus $\rho$ solves
a conservation law linear in $\rho$. Indeed the Newell-Daganzo
flow is linear in the free flow zone. Therefore, no shock is produced
inside the domain $[0,L]$ and characteristics are defined for all times.
In particular the value of $\rho$ is constant along  characteristics.
The characteristic exiting the domain at time $t$
enters the domain from the boundary at time $\tau(t)$.
Therefore we get $\rho(t,L)=\rho(0,\tau(t))=\frac{\Inflow(\tau(t))}{v(\tau(t))}$.
From \eqref{eq:Out} we get the desired conclusion.
\end{Proof}
\begin{remark}
This map is highly non-linear with respect to the control policy $v$
due to the definition of $\tau$. Hence, the classical techniques of linear control cannot be applied. Moreover, such formulation clearly shows 
how delays enter the input-output flow map.
The effect of the control $v$ at time $t$ on the outflow depends on the choice of $v$ on the time interval $[\tau(t),t]$, because of the presence
of the LET map in formula \eqref{eq:Outflow}.
\end{remark}
\subsection{Optimal control problem}\label{subsec:optimal_control}
We are now ready to define formally the problem of outflow tracking.
\begin{problem}\label{pb:optimal_control}
Let Hypothesis \ref{hyp:bounds} hold, fix $\ftarget \in \BV([0,T], [\fmin, \fmax])$ 
and $K>0$.
Find the control policy $v\in \BV([0,T], [\vmin,\vmax])$, with $\TV(v)\leq K$, which minimizes the functional $J: \BV([0,T], [\vmin,\vmax])\to\R$ defined by
\begin{equation}
	\label{eq:cost}
	J(v) := \int_{0}^{T}(\Outflow(t) - \ftarget(t))^2dt
\end{equation} 
where $\Outflow(t)$ is given by \eqref{eq:Outflow}.
\end{problem}
We prove later on, in Proposition \ref{prop:existence}, that Problem \ref{pb:optimal_control} admits a solution.
\begin{remark}
We use the same positive extreme values $\fmin,\; \fmax$ for both the inflow $\Inflow(.)$ and the target outflow $\ftarget(.)$ for simplicity of notation only. 
\end{remark}
\begin{remark}\label{rem:constant}
In the simple case where all the parameters are constant in time, i.e. $\Inflow, \; \Outflow, \; \ftarget, \; \rho_0$ do not depend on time, the problem has a a trivial solution which is $v =\dfrac{\ftarget}{\rho_0}$ realizing $J(v)=0$.
\end{remark}
\subsection{Cost variation as function of control policy variation}\label{subsec:gradient}
In this section we estimate the variation of the cost $J(v)$ with respect to the perturbations of the control policy $v$. This computation will allow to prove continuous dependence of the solution from the control policy.

We first fix the notation for integrals of $BV$ function with respect to Radon measures.
\begin{definition} \label{def:integrale}
Let $\phi$ be a $BV$-function and $\mu$ a Radon measure. We define
$$\int\phi(x^+)\,d\mu(x):=\int \phi(x)\,d\mu_c(x)+\sum_{i} m_i\phi(x_i^+),$$
where $\mu=\mu_c+\sum_{i} m_i \delta_{x_i}$ is the decomposition of $\mu$ into its continuous\footnote{We recall that any Radon measure on $\R$ can be decomposed into its continuous (AC+Cantor) and Dirac parts, as a consequence of the Lebesgue decomposition Theorem, see e.g. \cite{EGbook} .} and Dirac parts.
\end{definition}

We now compute the variation in the cost $J$ produced by needle-like variation
in the control policy $v(\cdot)$, i.e. variation of the value of $v(\cdot)$
on small intervals of the type $[t,t+\dt]$ in the same spirit as the needle variations
of Pontryagin Maximum Principle \cite{BPbook}.
\begin{definition}
Consider $v\in \BV([0,T], [\vmin, \vmax])$ and a time $t$ such that
$\tau^{-1}(0)=t_0\leq t<\tau(T)$ and $v(t^+)<\vmax$. 
Let $\dv>0$, $\dt>0$ be sufficiently small such that $t+\dt\leq \tau(T)$ and $v(t^+)+\dv\leq\vmax$.
We define a needle-like variation $v'(\cdot)$ of $v$, corresponding to $t$, $\dt$ and $\dv$
by setting $v'(s)=v(s)+\dv$ if $s\in [t,t+\dt]$ and $v'(s)=v(s)$ otherwise, see Figure \ref{fig:Needle_variation}.
\end{definition}
\begin{figure}[ht]
\centering
\begin{tikzpicture}
\draw [<->](4,0)-- (0,0) -- (0,4);
\draw(0,2).. controls (1,1) and (2,3) .. (3.5,2);
\draw(1,0)--(1,3);
\draw(2,0)--(2,3);
\draw(1,3)--(2,3);
\draw [dashed, thin](1,3) -- (0,3);
\node [below] (t) at (4,0) {$t$};
\node [left] (v) at (0,4) {$v$};
\node [below] (t1) at (1,0) {$t$};
\node [below] (dt) at (2,0.03) {$t+\Delta t$};
\node [left] (v1) at (0,2) {$v$};
\node [left] (dv) at (0,3) {$v'=v+\Delta v$};
\end{tikzpicture}
\caption{Needle-like variation of the velocity $v$.}
\label{fig:Needle_variation}
\end{figure}
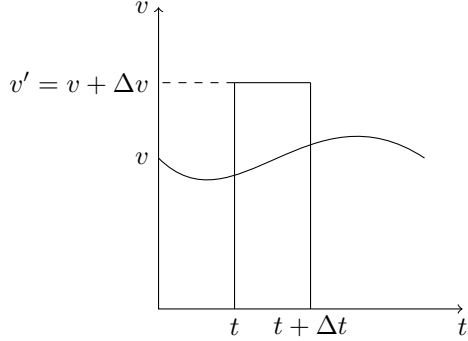
\begin{lemma}\label{prop:gradient}
 Consider $v\in \BV([0,T], [\vmin, \vmax])$ and let $v'$ be a needle-like variation of $v$.
 Then it holds:
\begin{equation}\label{eq:gradient}
\begin{split}
 \lim_{\dv\rightarrow 0^+}&\lim_{\dt\rightarrow 0^+}\dfrac{J(v')-J(v)}{\dv}=\\
& = 2\rho^2(t^-,L) v(t^+)-2\rho(t^-,L)\ftarget(t^+) + \\
& -\int_{0^+}^L v^2((t+s(x))^+)\,d\rho_{x}^2(t) + 2\int_{0^+}^L \ftarget((t+s(x))^+)v((t+s(x))^+)\,d\rho_{x}(t)+\\
& + 2 \dfrac{\Inflow(t^-)}{v(t^+)}\Big(\ftarget(t^+)-\frac{v(\tau^{-1}(t')^-)}{v(t^+)}\Inflow(t^-)\Big),
\end{split}
\end{equation}
where integrals are defined according to Definition \ref{def:integrale}.
For $\dv<0$, the limit for $\dv\to 0^-$ satisfies the same formula with
right limits replaced by left limits in the two integral terms in \eqref{eq:gradient}.
\end{lemma}
\begin{remark}\label{rem:gradient}
Notice that the condition $\tau^{-1}(0)=t_0<t$ implies that the outflow $\Outflow(s)$ $\in [t, t+\dt]$, depends only on the inflow $\Inflow(.)$ and not on the initial density $\rho_0.$ 
If such condition is not satisfied, the perturbation given by $\dv$ has a comparable effect on 
$\Outflow(.)$, but it needs to be estimated in two parts: one with respect to $\Inflow([0,t+\dt])$ and one with respect to $\rho_0(0,L-l)$ with $l$ being such that 
$$
\int_{0}^{t} v(s) ds = l.
$$
The condition $t+\dt \leq \tau(T)$ means that the perturbation $\dv$ has influence on the whole outflow $\Outflow(s)$ in the interval $[t, \tau^{-1}(t+\dt)]$. If this is not satisfied, then the influence of the perturbation is stopped at $T<\tau^{-1}(t+\dt) $, hence the variation $\Outflow(s)$ is smaller.
\end{remark}

\begin{Proof}
Let $\tau(t)$ be defined according to \eqref{eq:tau} and an outflow $\Outflow(t)$ according to \eqref{eq:Outflow}. 
For simplicity we assume that $v(\cdot)$ has a constant value $\vhat:=v(t^+)$
on $[t,t+\dt]$, the general case holding because of properties of $\BV$ functions.\\
We define $t'=t+\dt$ and $s'$ to be the unique value satisfying 
$$
\int_{0}^{s'} v(t'+\sigma)d\sigma = L - (\vhat+\dv)\dt,
$$
$s''$ to be the unique value satisfying
$$
\int_{0}^{s''} v(t'+\sigma)d\sigma = L - \vhat\dt,
$$
and $s'''=\tau^{-1}(t')-t'$, hence $\int_0^{s'''} v(t'+\sigma)d\sigma=L$. Notice that $s'<s''<s'''$. We also define the function
\begin{equation}
	\label{eq:x_s}
	x(s) = L-\int_{0}^s v(t'+\sigma)d\sigma.
\end{equation}
Remark that $x(s)$ is a decreasing function, with $x(0)=L$ , $x(s')=(\vhat+\dv)\dt$, $x(s'')=\vhat\dt$ and $x(s''')=0$.
\begin{figure}[ht]
\centering
\begin{tikzpicture}[scale=1.4,baseline=0]
\draw[thick](0,0)--(4,0);
\draw[thick](0,0)--(0,4.5);
\draw[thick](4,0)--(4,4.5);
\draw[dotted, thick](4,1.5)--(0,1.5);
\draw[<->](2,2.3)--(2,1.5);
\draw[<->](4.1,2.3)--(4.1,1.5);
\draw[<->](4.1,2.3)--(4.1,2.8);
\draw[<->](4.1,3.3)--(4.1,2.8);
\draw[<->](4.1,3.3)--(4.1,3.8);
\draw[dotted,thick](4,2.3)--(0,2.3);
\draw[dotted,thick](0,2.8)--(4,2.8);
\draw[dotted,thick](0,3.3)--(4,3.3);
\draw[dotted,thick](0,3.8)--(4,3.8);
\draw[dotted](2.5,0)--(2.5,3.3);
\draw[dotted](1.5,0)--(1.5,2.8);
\draw[white, fill=gray!20!white, opacity=0.5](0,0)rectangle(4,4.5);
\draw[line width=1pt] (0,0.2)--(1.5,0.7)--(3,1)--(4,1.5);
\draw[line width=1pt] (0,2.3)--(1.5,2.8)--(2.5,3.3)--(4,3.8);
\node[below]at(-0.1,0){$0$};
\node[below]at(4.1,-0){$L$};
\node[below]at(1.6,0){$L-x(s')$};
\node[below]at(2.8,0){$L-x(s'')$};
\node[left]at(0,0.2){$\tau(t)$};
\node[left]at(0,2.3){$t'=t+\Delta t$};
\node[left]at(0,2.8){$t'+s'$};
\node[left]at(0,3.3){$t'+s''$};
\node[left]at(0,3.8){$t'+s'''$};
\node[right]at(4,1.4){$t$};
\node[left]at(0,1.5){$t$};
\node[right]at(4,4){$\tau^{-1}(t+\Delta t)=\tau^{-1}(t')$};
\node[right]at(2,1.95){$\Delta t$};
\node[right]at(4.1,1.95){$I_1$};
\node[right]at(4.1,2.5){$I_2$};
\node[right]at(4.1,3){$I_3$};
\node[right]at(4.1,3.5){$I_4$};
\end{tikzpicture}
\caption{Graphical representation for the notation used in section \ref{subsec:gradient}}
\label{fig:variation}
\end{figure}
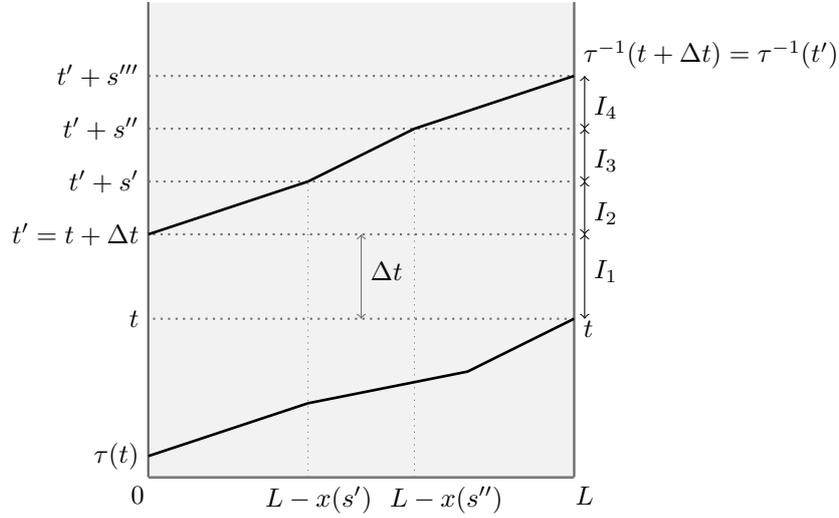
We denote with $\Outflow'(s)$ the outflow, $\tau'(s)$ the LET (see \eqref{eq:tau}) and $\rho'(s,x)$ the density for the policy $v'$. Clearly, we have $\Outflow'(s) = \Outflow(s)$ for $s \in [0,t] \cup [\tau^{-1}(t'),T]$ and $\tau'(s) = \tau(s)$ for $s\in [t_0, t] \cup [\tau^{-1}(t'),T]$.\\
To compute the variation, we distinguish four time intervals:
$I_1 = (t,t')$, $I_2 = (t',t'+s')$, $I_3 = (t'+s',t'+s'')$ and $I_4 = (t'+s'',\tau^{-1}(t')),$ see Figure \ref{fig:variation}.
The variation of the cost in the first interval can be directly computed
as function of the velocity variation, while in the other intervals
the delays in the outflow formula \eqref{eq:Outflow} will render
the computation more involved.
We denote with $J_1,\ldots,J_4$ the contributions to $\lim_{\dt\to 0^+}(J(v')-J(v))/\dv$ in the four intervals and estimate them separately.\\[0.2cm]
\textbf{CASE 1} : $I_1 = (t,t')$. 
Let $s\in [0,t'-t] =[0,\dt]$, then 
$\Outflow(t+s)=\rho(t, L-s\vhat)\vhat$ and $\Outflow'(t+s)=\rho(t, L-s(\vhat+\dv))(\vhat+\dv)$.
We have:
\begin{equation}\label{eq:J1}
J_1=\lim_{\dt\rightarrow 0^+} \dfrac{1}{\dt} \Big[ \int_{0}^{\dt}\Big (\Outflow'(t+s)-\ftarget(t+s)\Big)^2 ds - \int_{0}^{\dt}\Big (\Outflow(t+s)-\ftarget(t+s)\Big)^2 ds \Big] =
\end{equation}
$$\lim_{\dt\rightarrow 0^+} \dfrac{1}{\dt} \Big[ \int_{0}^{\dt}\Outflow'^2(t+s)-\Outflow^2(t+s)-2\ftarget(t+s)\Big(\Outflow'(t+s)-\Outflow(t+s)\Big)ds \Big] =$$
Substituting the expressions for the outflows we get
$$\lim_{\dt\rightarrow 0^+} \dfrac{1}{\dt} \Big[ \int_{0}^{\dt} \rho^2(t, L-s(\vhat+\dv))(\vhat+\dv)^2 - \rho^2(t, L-s\vhat)\vhat^2 ds + $$
$$ - \int_{0}^{\dt}2\ftarget(t+s)\Big( \rho(t, L-s(\vhat+\dv))(\vhat+\dv)-\rho(t, L-s\vhat)\vhat \Big)ds\Big]=$$
Dividing the first integral in two parts and making the change of variable
$ \sigma=s\dfrac{\vhat+\dv}{\vhat}$
$$\lim_{\dt\rightarrow 0^+} \dfrac{1}{\dt} \Big[ \int_{0}^{\dt(1+\frac{\dv}{\vhat})} \rho^2(t, L-\sigma\vhat)(\vhat+\dv)^{\cancel{2}}\dfrac{\vhat}{\cancel{\vhat+\dv}}d\sigma - \int_{0}^{\dt} \rho^2(t, L-s\vhat)\vhat^2 ds +$$
$$ -\int_{0}^{\dt}2\ftarget(t+s)\Big(\vhat (\rho(t, L-s(\vhat+\dv))-\rho(t, L-s\vhat))+\dv(\rho(t, L-s(\vhat+\dv)))\Big) ds\Big] =$$
After simple algebraic manipulation we get:
$$ \lim_{\dt\rightarrow 0^+} \dfrac{1}{\dt} \Big[ \int_{0}^{\dt(1+\frac{\dv}{\vhat})} \rho^2(t, L-s\vhat)\dv\vhat ds+ \int_{\dt}^{\dt(1+\frac{\dv}{\vhat})}\rho^2(t, L-s\vhat)\vhat^2 ds + $$
$$- \int_{0}^{\dt}2\ftarget(t+s)\Big(\vhat (\rho(t, L-s(\vhat+\dv))-\rho(L-s\vhat))+\dv(\rho(t, L-s(\vhat+\dv)))\Big) ds\Big] =$$
$$ \lim_{\dt\rightarrow 0^+} \dfrac{1}{\dt} \Big[ \int_{0}^{\dt}\rho^2(t, L-s\vhat)\dv\vhat ds + \int_{\dt}^{\dt(1+\frac{\dv}{\vhat})}\rho^2(t, L-s\vhat)(\vhat^2+\dv\vhat) ds $$
$$-  \int_{0}^{\dt}2\ftarget(t+s)\Big(\vhat (\rho(t, L-s(\vhat+\dv))-\rho(t, L-s\vhat))+\dv(\rho(t, L-s(\vhat+\dv)))\Big) ds\Big] =$$
Taking the limit as $\dt\to 0^+$, we get:
$$\rho^2(t^-,L)\vhat\dv+\rho^2(t^-,L)\cancel{\vhat}(\vhat+\dv)\dfrac{\dv}{\cancel{\vhat}}+$$
$$-2\ftarget(t^+)[\vhat(\cancel{\rho(t^-,L)}-\cancel{\rho(t^-,L)})]-2\ftarget(t^+)\dv\rho(t^-,L)=$$
$$ \rho^2(t^-,L)\vhat\dv+\rho^2(t^-,L)(\vhat+\dv)\dv-2\ftarget(t^+)\dv\rho(t^-,L),$$
hence $$ J_1=2\rho^2(t^-,L)\vhat + \rho^2(t^-,L)\dv-2\ftarget(t^+)\rho(t^-,L),$$
thus
$$
\lim_{\dv\to 0^+} J_1=2\rho^2(t^-,L) v(t^+) -2\ftarget(t^+)\rho(t^-,L).
$$

\vspace{0.2cm}
\textbf{CASE 2} : $I_2 = (t',t'+s')$.
If $s \in [0, s']$ then $\Outflow(t'+s)=\rho(t',x(s))v(t'+s)$ and $\Outflow'(t'+s)=\rho((t',x(s)-\dv\dt))v(t'+s)$. After decomposing $J_2$ as done for $J_1$ in \eqref{eq:J1} and plugging in the expression of the outflows, we have 
\begin{equation}\label{eq:J2}
\begin{split}
J_2=\lim_{\dt\rightarrow 0^+} &\dfrac{1}{\dt}\Big[ \int_{0}^{s'} v^2(t'+s)\Big(\rho^2(t',x(s)-\dv\dt)-\rho^2(t',x(s))\Big)ds +\\
& - \int_{0}^{s'}2\ftarget(t'+s)v(t'+s)\Big(\rho(t',x(s)-\dv\dt)-\rho(t',x(s))\Big)ds\Big].
\end{split}
\end{equation}
Applying the change of variable $s \rightarrow x(s)$ (see \eqref{eq:x_s}), it holds
$$ J_2= \lim_{\dt\rightarrow 0^+}\dfrac{1}{\dt}\Big[  \int_{0^+}^L v^2(t'+s(x))\Big(\rho^2(t',x-\dv\dt)-\rho^2(t',x)\Big)dx  +$$ $$-\int_{0^+}^{L}2\ftarget(t'+s(x))v(t'+s(x))\Big(\rho(t',x-\dv\dt)-\rho(t',x)\Big)dx\Big]. $$
Notice that this change of variable is justified by Lemma \ref{lem:lim_append} of the Appendix. 
Using Lemma \ref{lem:BVappend} of the Appendix, we get:
 $$\lim_{\dv\to 0^+} J_2=-\int_{0^+}^L v^2((t'+s(x))^+)\,d\rho_{x}^2(t',x)$$
 $$ + 2\int_{0^+}^L \ftarget((t'+s(x))^+)v((t'+s(x))^+)\,d\rho_{x}(t',x).$$
\vspace{0.2cm}
 \textbf{CASE 3} : $I_3 = (t'+s',t'+s'')$.
If $s \in [s',s'']$  then $\Outflow(t'+s)=\rho(t',x(s))v(t'+s)$ and 
$$\Outflow'(t'+s)=v(t'+s)\dfrac{g(s)}{\vhat+\dv}, \qquad
g(s)=\Inflow\left(t'-\dfrac{x(s)}{\vhat+\dv}\right).$$ 
After decomposing $J_3$ as done for $J_1$ in \eqref{eq:J1} and plugging in the expression of the outflows, we get
 $$\lim_{\dt\rightarrow 0^+} \dfrac{1}{\dt} \Big[ \int_{s'}^{s''}v^2(t'+s)\dfrac{g^2(s)}{(\vhat+\dv)^2} - \rho^2(t',x(s))v^2(t'+s)+$$
 $$-2\ftarget(t'+s)\Big(  v(t'+s)\dfrac{g(s)}{\vhat+\dv} - \rho(t',x(s))v(t'+s) \Big)\Big]ds=$$
 Observe that $\lim_{\dt\to 0^+}s'=\lim_{\dt\to 0^+}s''=\tau^{-1}(t')^--t'$ and $\int_{s'}^{s''}v(t'+\sigma)d\sigma=\dv\dt$, then
 \begin{equation}
 \begin{split}
\dv\ J_3=&\dfrac{\dv}{v(\tau^{-1}(t')^-)} v^2(\tau^{-1}(t')^-)\Inflow^2(t'^-)\Big[\Big(\dfrac{1}{\vhat+\dv}\Big)^2-\Big(\dfrac{1}{\vhat}\Big)^2\Big] - \\
&\dfrac{\dv}{v(\tau^{-1}(t')^-)} 2\ftarget(\tau^{-1}(t')^-)v(\tau^{-1}(t')^-)\Inflow(t'^-)\Big( \dfrac{1}{\vhat+\dv}-\dfrac{1}{\vhat}\Big)\Big], 
\end{split}
 \end{equation}
thus $$\lim\limits_{\dv\to 0^+}J_3=0.$$\vspace{0.2cm}
 \textbf{CASE 4} : $I_4 = (t'+s'',t'+s''')$. If $s \in [s'',s''']$  then we compute $$\Outflow(t'+s)=\dfrac{h(s)}{\vhat}v(t'+s) \qquad h(s) = \Inflow\left(t'-\dfrac{x(s)}{\vhat}\right)$$ and 
 $$\Outflow'(t'+s)=v(t'+s)\dfrac{g(s)}{\vhat+\dv} \qquad g(s)=\Inflow\left(t'-\dfrac{x(s)}{\vhat+\dv}\right).$$ 
We decompose $J_4$ as done with $J_1$ in \eqref{eq:J1}, plug in the expression of the outflows, and use the equality $\int_{s''}^{s'''}v(t'+\sigma)\,d\sigma=\hat v$. The, denoting $\tilde v=v(\tau^{-1}(t')^-)$, we have
 $$\dv J_4=\dfrac{\vhat}{\tilde v} \Big[ \tilde v^2\Inflow^2(t'^-)\Big[\Big(\dfrac{1}{\vhat+\dv}\Big)^2-\Big(\dfrac{1}{\vhat}\Big)^2\Big] - 2\ftarget(\tau^{-1}(t')^-)\tilde v \Inflow(t'-)\Big[ \dfrac{1}{\vhat+\dv}-\dfrac{1}{\vhat}\Big]\Big].$$
By passing to the limit, we get
 \begin{equation*}
\lim_{\dv\to 0^+}J_4=2\ftarget(\tau^{-1}(t')^-)\frac{\Inflow(t'-)}{\hat v}-2\frac{\tilde v}{\hat v^2} \Inflow(t'-)^2. 
 \end{equation*}
 \vspace{0.1cm}
\end{Proof}


Lemma \ref{prop:gradient} and Remark \ref{rem:gradient} allow us to prove the following:
\begin{proposition}\label{prop:grad}
For every $K>0$ and $C>0$,
the functional $J$ is Lipschitz continuous on 
$\Omega = \{ v \in\BV : \TV(v) \leq K \} \cap \{ v \in \Lspace{\infty}: \norma{v}_{\infty} \leq C \}$
endowed with the norm $\norma{v}_{\Lspace{1}}$.
\end{proposition}
\begin{Proof}
Let $v$, $\tilde{v}\in\Omega$. Then $v-v'$ is in $\BV$ and can be approximated by
piecewise constant functions. This means the $v-v'$ can be approximated in $\BV$ by needle-like
variations as in Lemma \ref{prop:gradient}. The right-hand side of \eqref{eq:gradient}
is uniformly bounded (since $v\in\Omega$ and $\rho\in\BV$ with uniformly bounded variation).
Therefore we conclude that $\modulo{J(v)-J(v')}\leq  C \norma{v-v'}_{\Lspace{1}}$ for some $C>0$.
\end{Proof}

This allows to prove the following existence result.
\begin{proposition}\label{prop:existence}
Problem \ref{pb:optimal_control} admits a solution.
\end{proposition}
\begin{Proof}
The space $\Omega = \{ v \in\BV : \TV(v) \leq K \} \cap \{ v \in \Lspace{\infty}: \norma{v}_{\infty} \leq C \}$ is compact in $\Lspace{1}$, see e.g. \cite{AFPbook}, and $J$
is Lipschitz continuous on $\Omega$, thus there exists a minimizer of Problem \ref{pb:optimal_control}.
%
\end{Proof}
\section{Control policies}\label{sec:policies}
In this section, we define three control policies for the time-dependent maximal speed $v$. The first, called the instantaneous policy (IP), is defined by minimizing the instantaneous contribution for the cost $J(v)$ at each time. We will show that such control policy does not provide a global minimizer, due 
to delays in the control effect on the cost for the Problem \ref{pb:optimal_control}.
In particular, due to the bound $v\in [\vmin, \vmax])$ the instantaneous minimization may induce a larger cost at subsequent times. Then, we introduce a second control policy, called random exploration (RE) policy.
Such policy uses a random path along a binary tree, which correspond the upper and
lower bounds for $v$, i.e. $v=\vmax$ and $v=\vmin$.
Lastly, a third control policy is based on a gradient method descent (GDM) and uses the expression given by Lemma \ref{prop:gradient} to numerically find the gradient of the cost.
\subsection{Instantaneous policy}\label{subsec:InstPolicy}
\begin{definition}\label{def:instPolicy}
Consider Problem \ref{pb:optimal_control}. Define the \textbf{instantaneous policy} as follows:
\begin{equation}
	\label{eq:Ip}
	v(t) := P_{[\vmin,\vmax]}\Big( \ftarget(t^-)\cdot \dfrac{v(\tau(t)^-)}{\Inflow(\tau(t)^-)}\Big),
\end{equation}
where the projection $ P_{[\vmin,\vmax]}: \R\rightarrow\R$ is the function
\begin{equation}
	\label{eq:projection}
	P_{[a,b]}(x) := \left\{ \begin{array}{ll}
	a \quad \text{for} \: x<a,\\
	x \quad \text{for} \: x\in [a,b],\\
	b \quad \text{for} \: x>b.
	\end{array}
	\right.
\end{equation}
\end{definition} 
Notice that this would be the optimal choice if $\ftarget$ and $In$ would be constant,
see Remark \ref{rem:constant}.
The instantaneous policy can also be written directly in terms of the input-output map defined in  Proposition \ref{prop:inputOutputMap}.
 As we will show later, the instantaneous policy is not optimal in general, i.e., it does not provide an optimal solution $v$ for Problem \ref{pb:optimal_control}. Clearly, it provides the solution in the case of $\vmin$ sufficiently small and $\vmax$ sufficiently big so that the projection operator reduces to the identity, i.e., $v(t) = P_{[\vmin,\vmax]} \Big( \dfrac{\ftarget (t^-)}{\rho(L^-)}\Big)=\dfrac{\ftarget (t^-)}{\rho(L^-)}$ for all times. Indeed, in this case the output $\Outflow(t)$ coincides with $\ftarget(t)$, hence the cost $J(v)$ is zero.
\subsection{Random exploration policy}
The random exploration policy is defined as follows:
\begin{definition}\label{def:forecastPolicy}
Given the extreme values for the maximal speed, $\vmax$ and $\vmin$, and a time step $\dt$,
the \textbf{random exploration policy} draws sequences of velocities from the set
$\{\vmax,\vmin\}$ corresponding to control policy values on the intervals $[i\dt,(i+1)\dt]$.
\end{definition}
Notice that maximal speeds according to this algorithm can be generated for all times,
independently of the corresponding solution, in contrast to the instantaneous policy
which is based on the maximal speed at previous times.
We will use numerical optimization to choose the best among the generated random policies, showing in particular that the instantaneous policy is not optimal in general. 

\subsection{Gradient method}
We use needle-like  variations and the analytical expression in \eqref{eq:gradient} to numerically 
compute one-sided variations of the cost. We consider such variations as estimates of the
gradient of the cost in $L^1$.
More precisely, we give the following definition.
\begin{definition}\label{def:gradientPolicy}
The \textbf{gradient policy} is the result of a first-order optimization algorithm to find a local minimum to  Problem \ref{pb:optimal_control} using the Gradient Descent Method and the expression in \eqref{eq:gradient}, 
stopping at a fixed precision tolerance.
\end{definition}
We will show that the gradient method gives very good results compared to the other policies
taking into account the computational complexity.
\section{Numerical simulations}\label{sec:num_simulations}
In this section we show the numerical results obtained by implementing the policies described in section \ref{sec:policies}. The numerical algorithm for all the approaches is composed of two steps:
\begin{enumerate}
\item Numerical scheme for the conservation law \eqref{eq:cons_law}. The density values are computed using the classical Godunov scheme, introduced in \cite{G59}.
\item Numerical solution for the optimal control problem, i.e., computation of the maximal speed using the instantaneous control, random exploration policy and gradient descent.
\end{enumerate}
Let $\dx$ and $\dt$ be the fixed space and time steps, and set $\xj{j+}{1} = j\dx$, the cell interfaces such that the computational cell is given by $C_j = [ \xj{j-}{1}, \xj{j+}{2}].$ The center of the cell is denoted by $x_j = (j-\dfrac{1}{2})\dx$ for $j\in \Z$ at each time step $t^n = n\dt$ for $n\in \N.$ We fix $\mathcal{J}$ the number of space points and $T$ the finite time horizon.
We now describe in detail the two steps.
\subsection{Godunov scheme for hyperbolic PDEs}
The Godunov scheme is a first order scheme, based on exact solution to Riemann problems. Given $\rho(t,x),$ the cell average of $\rho$ in the cell $C_j$ at time $t^n$ is defined as 
\begin{equation}
	\label{eq:cell_avg_god}
	\rho_j = \dfrac{1}{\dx} \int_{\xj{j-}{1}}^{\xj{j+}{1}} \rho(t^n,x)dx.
\end{equation}
Then, the Godunov scheme consists of two main steps:
\begin{enumerate}
\item Solve the Riemann problem at each cell interface $\xj{j+}{1}$ with initial data $(\rho_j, \rho_{j+1})$.
\item Compute the cell averages at time $t^{n+1}$ in each computational cell and obtain $\rho_j$.
\end{enumerate}
\begin{remark}
Waves in two neighboring cells do not intersect before $\dt$ if the following CFL (Courant-Friedrichs-Lewy) condition holds:
\begin{equation}
	\label{eq:CFL}
	\dt\max_{j\in\Z}{\modulo{f'(\rho_j)}}\leq \dfrac{1}{2}\min_{j\in\Z}\dx.
\end{equation}
\end{remark}
The Godunov scheme can be expressed in conservative form as:
\begin{equation}
	\label{eq:god_scheme}
	\rho_j^{n+1}= \rho_j^n - \dfrac{\dt}{\dx}\Big( F(\rho_j^n, \rho_{j+1}^n,v^n) -  F(\rho_{j-1}^n\rho_j^n,v^n)\Big)
\end{equation}
where $v^n$ is the maximal speed at time $t^n$. Additionally, $F(\rho_j^n, \rho_{j+1}^n,v^n)$ is the Godunov numerical flux that in general has the following expression:
\begin{equation}
	\label{eq:god_flux}
	F(\rho_j^n, \rho_{j+1}^n,v^n) = \left\{ \begin{array}{ll}
	\min_{z\in [\rho_j^n,\rho_{j+1}^n]} f(z,v^n)& \qquad \text{if } \rho_j^n\leq \rho_{j+1}^n,\\
	\max_{z\in \rho_{j+1}^n, \rho_j^n} f(z,v^n) &\qquad \text{if } \rho_{j+1}^n \leq \rho_j^n.
	\end{array}
	\right.
\end{equation}
For clarity, we included as an argument for the Godunov scheme the maximal velocity so that the dependence of the scheme on the optimal control could be explicit.
\subsection{Velocity policies}
The next step in the algorithm consists of computing a control policy $v$ that can be used in the Godunov scheme with the different approaches introduced in section \ref{sec:policies}. In particular, for the instantaneous policy approach we compute the velocity at each time step using the instantaneous outgoing flux. Instead, using the other two approaches,  the RE and the GDM, we compute beforehand the value of the velocity at each time step and then use it to solve the conservation law with the Godunov scheme.
\subsubsection{Instantaneous policy}
We follow the control policy described in section \ref{subsec:InstPolicy} for the instantaneous control. At each time step, the velocity $v^{n+1}$ is computed using the following formula:
\begin{equation}
	\label{eq:inst_velocity}
	v^{n+1}=v(t^{n+1}) = P_{[\vmin,\vmax]} \Big( \dfrac{\ftarget(t^{n})}{\rho_{\mathcal{J}}^{n}}\Big).
\end{equation}
\subsubsection{Random exploration policy}
 To compute for each time step the value of the velocity, we use a randomized path on a binary tree, see  Figure \ref{fig:bin_tree}. With such technique, we obtain several sequences of possible velocities. For each sequence the velocities are used to compute the fluxes for the numerical simulations. We then choose the sequence that minimizes the cost.
 \begin{figure}[ht]
	\centering
	\resizebox{0.4\columnwidth}{!}{
	\begin{tikzpicture}[level distance=1.5cm,
  level 1/.style={sibling distance=3cm},
  level 2/.style={sibling distance=1.5cm}]
  \node [blue]{$\vmax$}
    child {node [blue]{$\vmax$}
      child {node [blue]{$\vmax$}}
      child {node [red]{$\vmin$}}
    }
    child {node [red]{$\vmin$}
    child {node [blue]{$\vmax$}}
      child {node [red]{$\vmin$}}
    };
\end{tikzpicture}}
	\label{fig:bintree1}
	\resizebox{0.4\columnwidth}{!}{
	\begin{tikzpicture}[level distance=1.5cm,
  level 1/.style={sibling distance=3cm},
  level 2/.style={sibling distance=1.5cm}]
  \node [red]{$\vmin$}
    child {node [blue]{$\vmax$}
      child {node [blue]{$\vmax$}}
      child {node [red]{$\vmin$}}
    }
    child {node [red]{$\vmin$}
    child {node [blue]{$\vmax$}}
      child {node [red]{$\vmin$}}
    };
\end{tikzpicture}}
	\label{fig:bintree2}
	\caption{The first branches of the binary tree used for sampling the velocity.}
	\label{fig:bin_tree}
\end{figure}
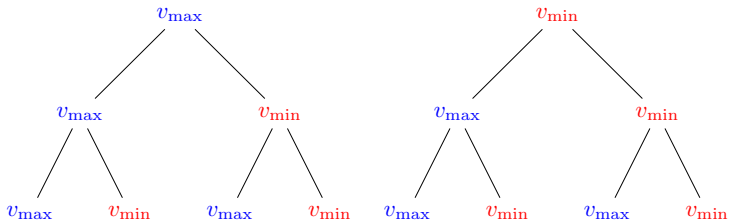
\begin{remark}
Notice that the control policy RE may have a very large variation, not respecting the bounds
given in Problem \ref{pb:optimal_control}. Therefore the found control policies may not be
feasible. 
However, we implement this technique for comparison with the results and performances obtained by the GDM.
\end{remark}
\subsubsection{Gradient descent method}
We first numerically compute one-sided variations of the cost using \eqref{eq:gradient}. Then, we use the classical gradient descent method \cite{Abook} to find the optimal control strategy and to compute the optimal velocity that fits the given outflow profile, as described in  Algorithm \ref{alg:gradientDescent}.%
\begin{algorithm}
\caption{\small{Algorithm for the gradient descent and computation of the optimal control}}\label{alg:gradientDescent}
\begin{algorithmic} \small{
	\STATE \textbf{Input data:} Initial and boundary condition for the PDE and initial velocity
	\STATE Fix a step tolerance $\epsilon$ and find a suitable step size $\alpha$ 
	\WHILE {$\modulo{J_{i+1}-J_i}\leq \epsilon$}
	\STATE Compute numerically cost variations $\nabla J_i $  
	\STATE Update the optimal velocity $v_{i+1}=v_i-\alpha\nabla J_i $
	\STATE Compute the new densities using Godunov scheme
	\STATE Compute the new value of the cost functional
	\ENDWHILE}%
\end{algorithmic}
\end{algorithm}
\subsection{Simulations}
We set the following parameters: $L=1, \; \mathcal{J}=100,\; T=15.0,\; \rho_{\mathrm{cr}}= 0.5,\; \rho_{\max}=1, \; \vmin = 0.5,\; \vmax= 1.0.$ Moreover, the input flux at the boundary of the domain is given by $\Inflow = \min{(0.3+0.3\sin(2\pi t^n),0.5)}$. We choose two different target fluxes $\ftarget = 0.3$ and $\ftarget = \modulo{(0.4\sin(t\pi - 0.3))}$. The initial condition is a constant density $\rho(0,x)=0.4.$ We use oscillating inflows to represent variations in typical inflow
of urban or highway networks at the 24h time scale.

\subsubsection{Test I: Constant Outflow}
\begin{figure}[ht]
	\centering
	\includegraphics[scale=0.3]{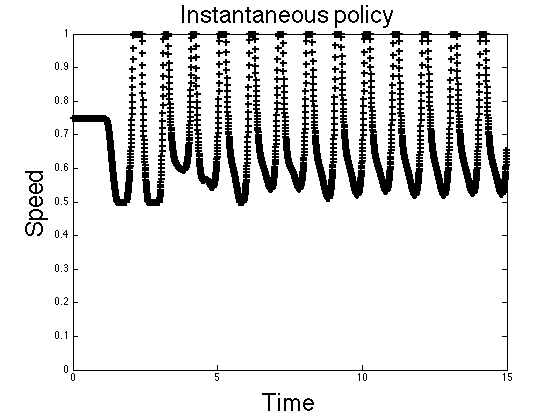}
	\includegraphics[scale=0.3]{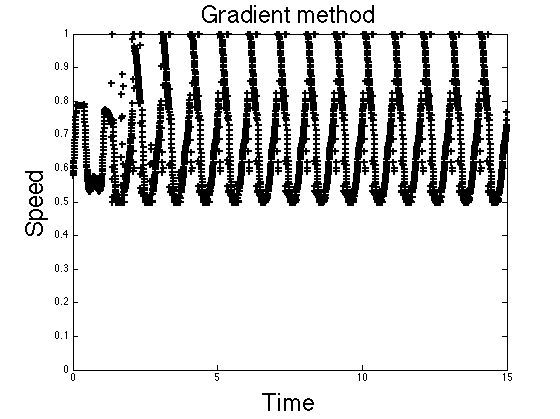}
	\caption{Speed obtained by using the instantaneous policy (left) and the gradient descent method (right) for a target flux $\ftarget=0.3$.}
		\label{fig:OptimalControl}
\end{figure}

\begin{figure}[ht]
	\centering
	\includegraphics[scale=0.4]{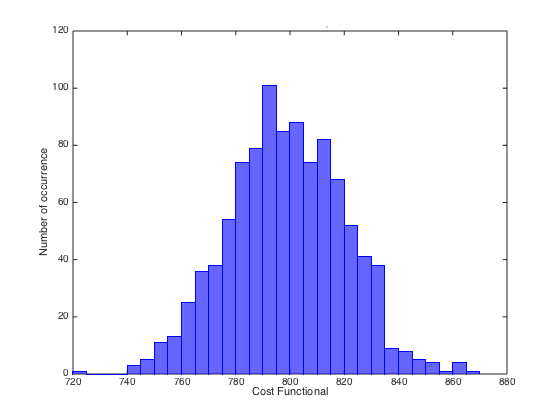}
	\caption{Histogram of the distribution of the value of the cost functional for the random exploration policy. We run 1000 different simulations.}
	\label{fig:histogram_const}
\end{figure}

\begin{figure}[ht]
\centering
\includegraphics[scale=0.3]{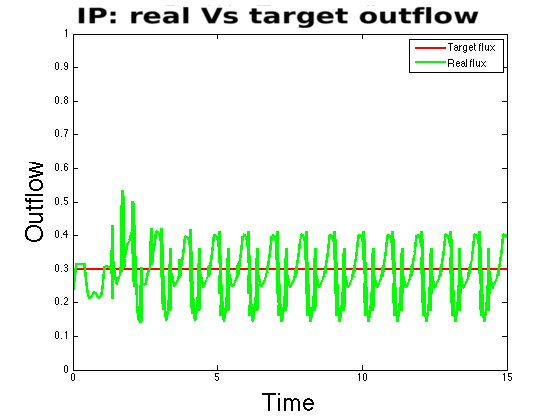}
\includegraphics[scale=0.3]{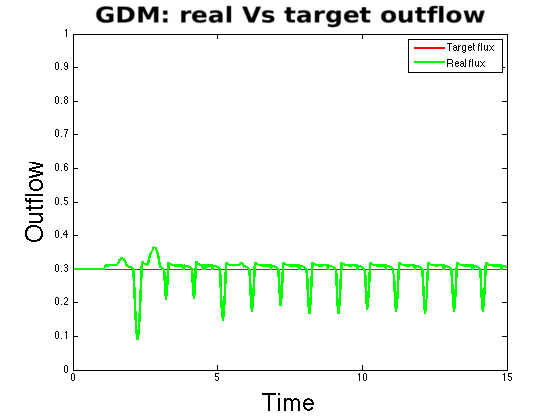}\\
\includegraphics[scale=0.3]{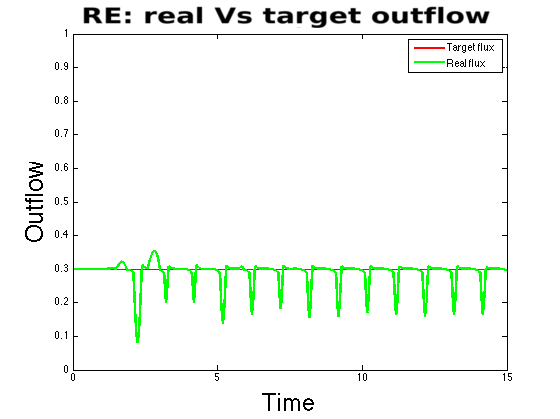}
\caption{Difference between the real outgoing flux and the target constant flux, computed with the instantaneous policy (top, left), the gradient method (top, right) and the random exploration policy (bottom).}
\label{fig:OutflowConst}
\end{figure}

In Figure \ref{fig:OptimalControl}, we show the time-varying speed obtained by using the instantaneous policy (left) and by using the gradient descent method (right). In each case, we notice that due to the oscillating input signal the control policy is also oscillating. We are aware, however, that from a practical point of view, the solution where the speed changes at each time step might be unfeasible. Nonetheless, these policies can be seen as periodic change of maximal speed for different time frames during the day when the time horizon is scaled to the day length. 
\begin{table}[ht]
\centering
\begin{tabular}{|c|c|}
\hline
Method & Cost Functional\\
\hline
Fixed speed $v=\vmax=1.0$ & $873.0786$\\
\hline
Fixed speed $v=\vmin=0.5$ & $785.2736$\\
\hline
Instantaneous policy & $850.3704$\\
\hline
Minimum of random exploration policy & $723.6733$\\
\hline
Gradient method& $735.0565$\\
\hline
\end{tabular}
\caption{Value of the cost functional for the different policies.}
\label{tab:costfunct}
\end{table}
In Table \ref{tab:costfunct}, we see the different results obtained for the cost functional computed at the final time for the different policies. For comparison, we also put the results of the simulations with a constant speed equal to the minimum and maximal velocity bounds. The instantaneous policy is outperformed by the random exploration policy and by the gradient method. For the random exploration policy, in the table we put the minimal value of the cost functional computed by the algorithm. In Figure \ref{fig:histogram_const} we can see the distribution of the different values of the cost functional over 1000 simulations. Moreover, in Figure \ref{fig:OutflowConst}, we can see the differences between the actual outflow obtained and the target one for all methods. 
We also compared the CPU time for the different simulations approaches (see Table \ref{tab:CPUtime}). As expected, the random exploration policy is the least performing while the instantaneous policy is the fastest one.
In addition, we computed the $\TV(v)$ for each one of the policies obtaining the following results:
\begin{itemize}
\item IP: $\TV(v)=12.6904$
\item RE: $\TV(v)=753.5$
\item GDM: $\TV(v)=70.81333.$
\end{itemize} 
\begin{table}[ht]
\centering
\begin{tabular}{|c|c|}
\hline
Method & CPU Time (s)\\
\hline
Instantaneous policy & $32.756$\\
\hline
Random exploration policy & $7577.390$\\
\hline
Gradient method& $1034.567$\\
\hline
\end{tabular}
\caption{CPU Time for the simulations performed with the different approaches.}
\label{tab:CPUtime}
\end{table}
\subsubsection{Test II: Sinusoidal Outflow}
\begin{figure}[ht]
	\centering
	\includegraphics[scale=0.3]{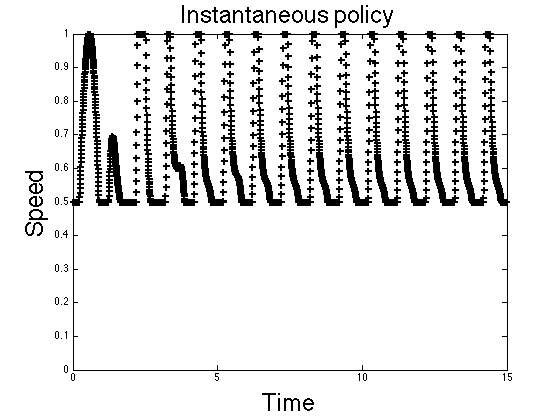}
	\includegraphics[scale=0.3]{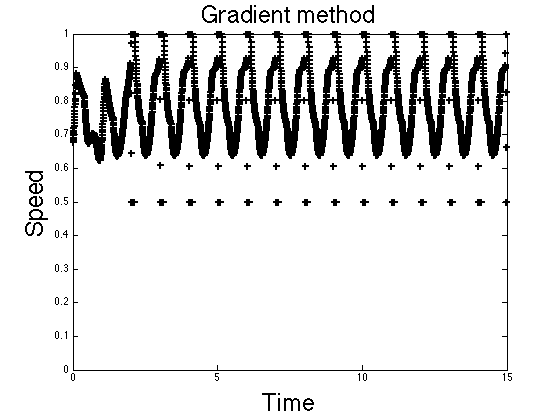}
	\caption{Speed obtained by using the instantaneous policy (left) and the gradient descent method (right) for a sinusoidal target flux.}
		\label{fig:SinusoidalOptimalControl}
\end{figure}

\begin{figure}[ht]
	\centering
	\includegraphics[scale=0.4]{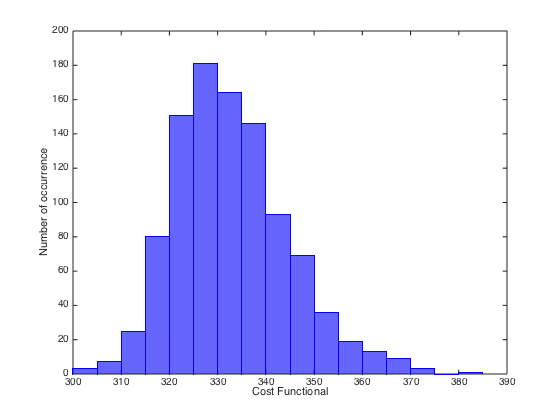}
	\caption{Histogram of the distribution of the value of the cost functional for the random exploration policy. We run 1000 different simulations.}
	\label{fig:histogram_sin}
\end{figure}

\begin{figure}[ht]
\centering
\includegraphics[scale=0.3]{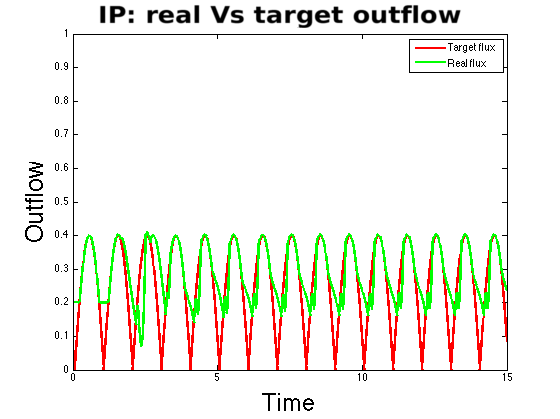}
\includegraphics[scale=0.3]{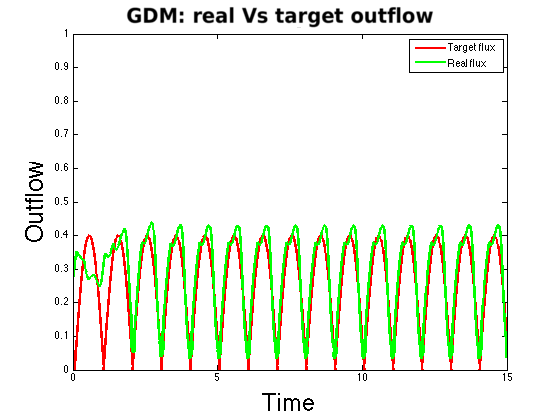}\\
\includegraphics[scale=0.3]{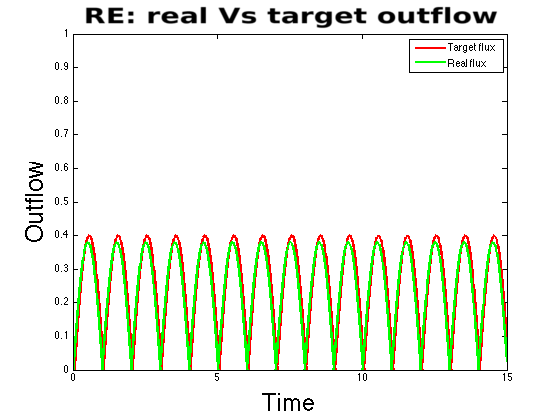}
\caption{Difference between the real outgoing flux and the target sinusoidal flux, computed with the instantaneous policy (top, left), the gradient method (top,right) and the random exploration policy (bottom).}
\label{fig:OutflowSin}
\end{figure}

In Figure \ref{fig:SinusoidalOptimalControl}, we show the optimal velocity obtained by using the instantaneous policy and by using the gradient descent method with a sinusoidal outflow. We show in Figure \ref{fig:histogram_sin}  the histogram of the cost functional obtained for the random exploration policy and in Figure \ref{fig:OutflowSin} we compare the real outgoing flux with the target one.
In Table \ref{tab:costfunct_II}, different results obtained for the cost functional computed at final time for the different policies are shown. Also in this case the instantaneous policy is outperformed by the other two. The CPU times give results similar to the previous test.
\begin{table}[ht]
\centering
\begin{tabular}{|c|c|}
\hline
Method & Cost Functional\\
\hline
Fixed speed $v=\vmax=1.0$ & $1.3979e+03$\\
\hline
Fixed speed $v=\vmin=0.5$ & $843.3395$\\
\hline
Instantaneous policy & $458.8874$\\
\hline
Minimum of random exploration policy & $303.8327$\\
\hline
Gradient method& $307.6889$\\
\hline
\end{tabular}
\caption{Value of the cost functional for the different policies.}
\label{tab:costfunct_II}
\end{table}
\section{Conclusions}\label{sec:conclusions}
In this work, we studied an optimal control problem for traffic regulation on a single road via variable speed limit. The traffic flow is described by the LWR model equipped with the Newell-Daganzo flux function. The optimal control problem consists in tracking a given target outflow in free flow conditions. We proved tje existence of a solution for the optimal control problem and provided explicit
analytical formulas for cost variations corresponding to
needle-like variations of the control policy. 
We proposed three different control policies design: instantaneous depending
only on the instantaneous downstream density, random simulations and gradient descent.
The latter, based on numerical simulations for the cost variation, represents the best
compromise between performance, computational cost and total variation of the control policy.\\
Future works will include the study of this problem in case of congestion and the extension to second order traffic flow models.
\section*{Appendix}
\begin{lemma}
\label{lem:lim_append}
Let $\alpha,\beta >0$, 
$\varphi, \psi \in \BV([0,L],\R)$, $\varphi>0$, define $x(s)=L-\int_0^s \varphi(\sigma)d\sigma$
and $\bar t$ by $\int_0^{\bar t} \varphi(\sigma)d\sigma=L- \alpha\dt$.
Using the change of variable $s\rightarrow x(s)$ we get:
\begin{equation} \label{eq:lim_prove}
\begin{split}
& \lim_{\dt\rightarrow 0^+} \dfrac{1}{\dt}\Big[ \int_{0}^{\bar t} \varphi^2(s)\Big(\psi (x(s)-\beta\dt)-\psi(x(s))\Big)ds \Big]=\\
 & \lim_{\dt\rightarrow 0^+}\dfrac{1}{\dt}\Big[  \int_{0^+}^L \varphi^2(s(x))\Big(\psi(x-\beta\dt)-\psi(x)\Big)dx \Big].
\end{split}
\end{equation}
\end{lemma}
\begin{Proof}
The change of variable $s \rightarrow x(s)$ inside the integral gives 
\begin{equation}
\begin{split}
&\lim_{\dt\rightarrow 0^+}\dfrac{1}{\dt}\Big[ \int_{0}^{\bar t} \varphi^2(s)\Big(\psi(x(s)-\beta\dt)-\psi(x(s))\Big)ds = \\
&\lim_{\dt\rightarrow 0^+}- \dfrac{1}{\dt} \int_{L}^{\alpha\dt} \varphi(s(x))\Big(\psi(x-\beta\dt)-\psi(x)\Big)dx =
\end{split}
\end{equation}
\begin{equation}\label{eq:limit_split}
\begin{split}
&\lim_{\dt\rightarrow 0^+}\dfrac{1}{\dt}\int_{0^+}^{L} \varphi(s(x))\Big(\psi(x-\beta\dt)-\psi(x)\Big)dx-\\
&\lim_{\dt\rightarrow 0^+}\dfrac{1}{\dt}\int_{0^+}^{\alpha\dt} \varphi(s(x))\Big(\psi(x(s)-\beta\dt)-\psi(x(s))\Big)dx.
\end{split}
\end{equation}
Observe that we need to specify the $0^+$ extremum in the integral, since the limit will provide Dirac terms inside the integral. We want now prove that the last addendum tends to zero. 
Denote by $\psi_x$ the distributional derivative of $\psi$, which is a measure, and decompose it as
in the continuous (AC+ Cantor) and Dirac part. By integrating $\psi_x$, we write $\psi=\tilde\psi+\sum_i m_i\chi_{[x_i,L]}$, with $\tilde\psi$ a continuous function, $m_i> 0$, $\sum_i m_i<+\infty$ and $x_i\in [0,L]$ . Hence, by the mean value theorem applied to $\tilde\psi$, we have 
\begin{equation}
\begin{split}
&\lim_{\dt\rightarrow 0^+}\dfrac{1}{\dt}\int_{0^+}^{\alpha\dt} \varphi(s(x))\modulo{\tilde \psi(x(s)-\beta\dt)-\tilde\psi(x(s))}dx\leq\\
&\lim_{\dt\rightarrow 0^+} \|\varphi\|_{\infty} \alpha \modulo{\tilde\psi(\tilde x -\beta\dt)-\tilde \psi(\tilde x)}=0,
\end{split}
\end{equation}
where $\tilde x\in (0,\alpha\dt)$ is a point (depending on $\dt$) and the limit is zero as a consequence of the continuity of $\tilde \psi$. The remaining term in \eqref{eq:limit_split} is then 
$$
\lim_{\dt\rightarrow 0^+}\dfrac{1}{\dt}\int_{0^+}^{\alpha\dt} \varphi(s(x))\sum_{x_i\in (0,\alpha\dt]} m_i (\chi_{[x_i-\beta\dt,L]}-\chi_{[x_i,L]})\,dx=$$
$$
\lim_{\dt\to 0^+} \frac{1}{\dt} \sum_{x_i\in (0,\alpha\dt]} \varphi(s(x_i)^-) m_i\beta\dt\leq \lim_{\dt\to 0^+} \beta \|\varphi\|_{\infty}  \sum_{x_i\in (0,\alpha\dt]}  m_i.
$$
Since $\psi$ is in $\BV$ the quantity $\sum_{x_i\in (0,\alpha\dt]}  m_i$ tends to zero
as $\dt$ tends to zero, thus we conclude.
\end{Proof}

\begin{lemma}
\label{lem:BVappend}
Let $\varphi, \psi \in \BV([a-\varepsilon,b+\varepsilon],\R)$, then 
\begin{eqnarray}\label{eq:app}
\lim_{\dt\rightarrow 0+} \dfrac{1}{\dt} \int_{a}^{b} \varphi(x)\Big( \psi(x-C \dt) -\psi(x)\Big)  dx =-C \int_a^b \varphi(x^+)d\psi_x(x),
\end{eqnarray}
where the integral in the right hand side is defined in Definition \ref{def:integrale}.
\end{lemma}
\begin{Proof}
We decompose the measure $\psi_x$ as $\psi_x=\ell\,d\lambda+\sum_i m_i\delta_{x_i}$,
where $\lambda$ is the Lebesgue measure, $\ell$ the Radon-Nikodym derivative of $\psi_x$
w.r.t. $\lambda$, $m_i>0$ and $\sum_i m_i<+\infty$. We approximate $\psi$ by piecewise
continuous functions $\psi^n$ defined as the integrals of $\psi^n_x=\ell\,d\lambda+\sum_{i\leq N(n)} m_i\delta_{x_i}$, where $N(n)$ is chosen such that $\sum_{i> N(n)} m_i<\frac{1}{n}$.\\
Define $I(n)=\cup_{i=1}^{N(n)}[x_i , x_i+C\dt] $ and by $I_c$ its complement in $[a,b]$.
Notice that for $x\in  [x_i , x_i+C\dt]$ we have $ \psi^n(x-C\dt) - \psi^n(x) = - m_i - \int_{x-C\dt}^{x}\ell\ d\lambda$ while on $I_c$ there are no jumps so 
$\psi^n(x-C\dt)-\psi^n(x) = - \int_{x-C\dt}^{x}\ell\ d\lambda$. We thus can write:
$$
\lim_{\dt\rightarrow 0+} \dfrac{1}{\dt} \int_{a}^{b} \varphi(x)\Big( \psi^n(x-C \dt) -\psi^n(x)\Big)  dx=
$$
\begin{equation*}
\begin{split}
\lim_{\dt\rightarrow 0^+}&\dfrac{1}{\dt}\sum\limits_{i=1}^{N(n)}\int_{x_i}^{x_i+C\dt}\varphi (x)\Big( \psi^n(x-C\dt)-\psi^n(x)\Big)dx +\\
&+ \int_{I_c}\varphi(x)\Big(\psi^n(x-C\dt)-\psi^n(x)\Big)dx=
\end{split}
\end{equation*}
\begin{equation}
\label{eq:int}
= \lim_{\dt\rightarrow 0^+}\dfrac{1}{\dt}\sum\limits_{i=1}^{N(n)}(-m_i) \int_{x_i}^{x_i+C\dt} \varphi(x)dx - \dfrac{1}{\dt}\int_{a}^{b} \varphi(x) \int_{x-C\dt}^{x}\ell\ d\lambda \, dx.
\end{equation}
Since $\varphi$ is in $\BV$ we can write:
$$
\lim_{\dt\rightarrow 0+} \dfrac{1}{\dt} \int_{a}^{b} \varphi(x)\Big( \psi^n(x-C \dt) -\psi^n(x)\Big)  dx= - \sum\limits_{i=1}^{N(n)} m_i \varphi (x^+)-\int_a^{b}\varphi(x)d(\ell\lambda) 
$$
$$
= - \int_a^{b} \varphi(x^+)d\Big(\sum_{i=1}^{N(n)} m_i \delta_{x_i} +\ell\lambda \Big) =- \int_a^{b} \varphi(x^+)d\psi^n_x$$
Now, the following estimates hold:
$$\modulo{\dfrac{1}{\dt}\int_{a}^{b} \varphi(x) \Big( \psi^n(x-C\dt)-\psi^n(x)\Big)dx - \dfrac{1}{\dt}\int_a^{b} \varphi(x)\Big( \psi(x-C\dt)-\psi(x)\Big)dx}$$
$$=\modulo{\dfrac{1}{\dt}\int_{a}^{b}\varphi \Big( \psi^n(x-C\dt) - \psi(x-C\dt) \Big) - \Big( \psi^n(x) - \psi(x)\Big)dx}$$
We can write $\psi^n(x-C\dt) = \psi(a) + \int_a^{x-C\dt}d\psi^n_x$ and $\psi(x-C\dt) = \psi(a) + \int_{a}^{x-C\dt}d\psi_x$, which gives us
$$ = \modulo{\dfrac{1}{\dt}\int_a^{b} \varphi (x) \Big( \int_a^{x-C\dt}dr_n - \int_{a}^{x} dr_n \Big)dx},$$
where $r_n=\psi-\psi^n$.
Taking the limit for $\dt\rightarrow 0^+$:
$$ \modulo{\dfrac{1}{\dt}\int_{a}^{b}\varphi(x)\Big( \psi^n(x-C\dt) - \psi^n(x)\Big)dx - \dfrac{1}{\dt}\int_a^{b} \varphi(x)\Big(\psi(x-C\dt) - \psi(x) \Big)dx}$$
$$\leq \modulo{\dfrac{1}{\dt}\int_a^{b} \varphi (x) \Big( -\int_{x-C\dt}^x dr_n\Big)dx} \leq$$
$$ \norma{\varphi}_{\infty} \dfrac{1}{\dt} \modulo{\int_a^{b}\int_{x-C\dt}^{x}dr_n dx} \leq \norma{\varphi}_{\infty} \dfrac{1}{n}.$$
The last inequality holds true because $\int_{x-C\dt}^{x}dr_n=\sum_i m_i \int_{x-C\dt}^{x}\ d\delta_{x_i}=$ \\$
\sum_i m_i \chi_{[x_i,x_i+C\dt]}$. Thus we get:
$$\lim_{\dt\rightarrow 0^+} \dfrac{1}{\dt} \int_a^{b} \varphi \Big(\psi(x-C\dt)-\psi(x)dx \Big)= \mathcal{O}\Big(\dfrac{1}{n}\Big) + \int_a^{b} \varphi (x^+)d\psi^n_x$$.
Let us now estimate the quantity
$$\modulo{\int_a^{b}\varphi(x^+)d\psi^n_x - \int_a^{b}\varphi(x^+)d\psi_x}.$$
Recalling that $\psi^n(x-C\dt) = \psi(a) + \int_a^{x-C\dt}d\psi^n_x$ and $\psi(x-C\dt) = \psi(a) + \int_{a}^{x-C\dt}d\psi_x$ we get 
$$\modulo{\int_a^{b} \varphi(x^+)d\Big(\sum\limits_{i\geq N(n)} m_i \delta_{x_i}\Big)}\leq \norma{\varphi}_{\infty}\dfrac{1}{n}.$$
Passing to the limit in $n$ we conclude.
\end{Proof}
\bibliographystyle{plain}
\bibliography{bibliotraffic}
\end{document}